\documentstyle[12pt]{article}
\title{Indexed identity and fuzzy set theory}
\author{Adonai S. Sant'Anna}
\date{{\it Dep. Matem\'atica, UFPR, C.P. 019081}\\{\it Curitiba, PR, 81531-990, Brazil}}
\begin{document}
\maketitle
%\tableofcontents
\newtheorem{definicao}{Definition}
\newtheorem{teorema}{Theorem}
\newtheorem{lema}{Lemma}
\newtheorem{exemplo}{Example}
\newtheorem{corolario}{Corollary}
\newtheorem{proposicao}{Proposition}
\newtheorem{axioma}{Axiom}
\newtheorem{observacao}{Remark}

\begin{abstract}
We introduce the concept of indexed identity, where the usual notion of identity is a particular case. Our mathematical framework allows us a generalized method for `indexing' predicates, which corresponds to `fuzzification' of properties, in an intuitive sense. It can be established a relationship between this indexed mathematics and fuzzy mathematics. We show some few examples that illustrate our ideas.
\end{abstract}

\section{Introduction}

	Fuzzy set theory appeared for the first time in 1965, in a famous paper by L. A. Zadeh \cite{Zadeh-65}. Since then a lot of fuzzy mathematics has been created and developed. Nevertheless, concepts like fuzzy set, fuzzy subset, and fuzzy equality (between two fuzzy sets) usually depend on the concept of {\em grade of membership\/}. This procedure does not allow us to define fuzzy equality between, e.g., two {\em Urelemente\/}, or even between two potatoes, up to the case that we consider a potato as a set.

	We present a concept of identity index which allows us to say how similar are two objects, even in the case that these two objects are not sets or fuzzy sets. Such a procedure allows us to define a sort of fuzzy membership (or grade of membership), fuzzy sets, fuzzy inclusion and fuzzy operations among fuzzy sets like union and intersection.

	This paper is the first one of a series dedicated to the concept of indexed identity and its applications.

\section{Indexed Identity}

	We use standard logical notation: $\neg$ is negation, $\wedge$ is conjunction, $\vee$ is disjunction, $\Rightarrow$ is conditional, and $\Leftrightarrow$ is biconditional. `$\exists$' and `$\forall$' denote the existence and universal quantifiers, respectively. Next, we define an indexed variables system by means of a {\em set-theoretical\/} predicate, following Suppes ideas about axiomatization \cite{Suppes-67}. By set theory we mean Zermelo-Fraenkel set theory (with or without {\em Urelemente\/}). But it is obvious that our system may be defined into the scope of other set theories.

\begin{definicao}\label{fvs}
An {\em indexed variables system\/} is an ordered pair $\langle X, \Xi\rangle$ that satisfies the following seven axioms:

\begin{description}

\item[F1] $X$ is a non-empty set.

\item[F2] $\Xi = \{\equiv_r\}_{r\in R}$ is a family of binary predicates defined on the elements of $X$, where $R$ is a subset of the interval of real numbers $[0,1]$, such that $1\in R$. When the ordered pair $(x,y)\in X\times X$ satisfies the binary predicate $\equiv_r$, we denote it by $x\equiv_r y$.

\item[F3] If $x\equiv_r y$ then $y\equiv_r x$.

\item[F4] $x\equiv_1 y$ iff $x=y$.

\item[F5] If $x\equiv_r y$ and $r\neq s$ then $\neg(x\equiv_s y)$.

\item[F6] $\forall x\forall y\exists r (x\equiv_r y)$.

\end{description}

\end{definicao}

\begin{definicao}\label{diff}
The {\em distinction\/} between two elements of $X$ is given by $D(x,y) = 1-r$ iff $x\equiv_r y$.
\end{definicao}

\begin{description}

\item[F7] $D(x,y)+D(y,z)\geq D(x,z)$

\end{description}

	For the sake of simplicity we can call $X$ as an indexed variables system or an indexed system. We call the binary relation $\equiv_r$ indexed identity with index $r$ or simply an indexed identity.

	Here follows an intuitive interpretation of the axioms and primitive concepts. $X$ is our space of variables. The sentence $x\equiv_r y$ corresponds to say that $x$ and $y$ do have an {\em identity index\/} $r$, where $0\leq r\leq 1$. If $r$ is 1 (one) then $x$ and $y$ are identical objects. If $r$ is not 1 then $x$ and $y$ are different objects. Nevertheless, if $r$ is a number `close' to 1, then $x$ and $y$ are very `similar' objects, i.e., `almost identical' objects. On the other hand, if $r$ is close to 0 (zero) then $x$ and $y$ are very different objects. They have `almost zero-identity'. This is reflected in the last axiom. According to {\bf F7} if $x$ and $y$ are very similar (index $r$ equal to 0.9, for example) and $y$ and $z$ are also very similar (index $r$ equal to 0.9, as another example), then $D(x,y) = D(y,z) = 0.1$. So, the distinction between $x$ and $z$ should be less or equal to 0.2, i.e., they should have an identity index greater or equal to 0.8. In the particular case where $D(x,y) = D(y,z) = 0$, then $D(x,z) = 0$, i.e., if $x=y$ and $y=z$ then $x=z$. Hence, axiom {\bf F7} is a generalization of the transitivity property of usual equality. The condition, in axiom {\bf F2} that $1\in R$ ensures consistency, since $x\equiv_1 y \Leftrightarrow x=y$.

	It is very important to emphasize that our mathematical framework is based on usual set theory, like Zermelo-Fraenkel's, for example.

\begin{teorema}
$D(x,y)$ is a distance between two points, and $\langle X, \Xi\rangle$ induces a metric space $\langle X, D\rangle$.
\end{teorema}

\noindent
{\bf Proof:} According to {\bf F4} and definition (\ref{diff}), we have $D(x,x) = 0$. According to {\bf F4} $D(x,y)>0$ if $x\neq y$. According to {\bf F3} $D(x,y) = D(y,x)$. According to {\bf F7} $D(x,z)\leq D(x,y)+D(y,z)$. So, $D(x,y)$ is the distance between $x$ and $y$; and $\langle X,D\rangle$ is a metric space.$\Box$\\

	Note that $D$ is a metric such that $D(x,y)\in[0,1]$.

\begin{teorema}
Any metric space $\langle M,d\rangle$ induces an indexed variables system $\langle M,\Xi\rangle$
\end{teorema}

\noindent
{\bf Proof:} If $\langle M,d\rangle$ is a metric space, then $d:M\times M\to \Re^+$ is a metric, where $\Re^+$ is the set of nonnegative real numbers. If we define a function $f:\Re^+\to [0,1]$ such that $f(x) = \frac{x}{1+x}$, then it is easy to prove that the function $D:M\times M\to [0,1]$ given by $D(a,b) = f(d(a,b))$ is a metric whose images belong to the set $[0,1)$. That is, we can define $D(a,b)$ as a distinction between $a$ and $b$ in the sense that $a\equiv_r b$ iff $D(a,b) = 1-r$, where $r\in (0,1]$. The reader can easily verify that $\langle M, \Xi\rangle$ is an indexed variables system, where $\Xi = \{\equiv_r\}_r$. $\Box$

\section{How to Calculate the Identity Index?}

	One natural question is: how to calculate the index $r$? In other words, what are our criteria to say how similar $x$ and $y$ are? One possible answer is the use of a family $\{A_i\}$ of unary predicates defined over the elements of $X$. For practical purposes we could consider a finite family $\{A_i\}_{i\in F}$, where $F$ is a finite set with cardinality, say, 100. If $x$ and $y$ are objects that share 93 predicates of the family $\{A_i\}_{i\in F}$, then they have 93\% of similarity, i.e., $x\equiv_{0.93} y$. If $x$ and $y$ have nothing in common, then they are totally different objects, i.e., $x\equiv_{0.00}y$. On the other hand, if $x$ and $y$ share all the predicates of the family $\{A_i\}_{i\in F}$, then $x$ and $y$ are identical objects, that is, they are the very same object with two different names or labels. Of course, this procedure just works if we make an adequate choice of possible predicates that objects of a given universe $X$ may (or not) satisfy. Such an adequate choice of predicates depends on the problem that we want to solve.

	It is important to say what do we mean by `two objects that share 93 predicates'. We say that $x$ and $y$ share one given predicate $A_i$ iff we have $A_i(x)\wedge A_i(y)$ or $\neg A_i(x)\wedge \neg A_i(y)$. We say that $x$ and $y$ do not share the predicate $A_i$ iff we have $A_i(x)\wedge \neg A_i(y)$ or $\neg A_i(x)\wedge A_i(y)$. We say that $x$ and $y$ share $n$ predicates iff there are $n$ different predicates $A_i$ that $x$ and $y$ share.

	As a simple example consider the set $X = \{2,3,8\}$, and the following family of predicates $\{A_1,A_2,A_3\}$, where $A_1(x)$ means that $x$ is an even number, $A_2(x)$ says that $x$ is an odd number and $A_3(x)$ corresponds to say that $x$ is a prime number. In this case we can calculate the index $r$ of indistinguishability between two elements $x$ and $y$ of $X$ as it follows:

\begin{equation}\label{calculandor}
r = \frac{\mbox{number of predicates $A_i$ that $x$ and $y$ share}}{\mbox{total number of predicates}}
\end{equation}

	So, we have $2\equiv_{1/3} 3$, $3\equiv_{0} 8$, $2\equiv_{2/3}8$, $2\equiv_1 2$, $3\equiv_1 3$, and $8\equiv_1 8$, since $A_1(2)\wedge \neg A_2(2)\wedge A_3(2)\wedge \neg A_1(3)\wedge A_2(3)\wedge A_3(3)\wedge A_1(8)\wedge \neg A_2(8)\wedge \neg A_3(8)$. Thus $D(2,3) = \frac{2}{3}$, $D(3,8) = 1$, and $D(2,8) = \frac{1}{3}$. It is easy to verify that $X$ is an indexed variables system.

	So, in this context, the number 2 is more similar to number 8, than to number 3. In this same context, numbers 3 and 8 have nothing in common, since their index of indistinguishability $r$ is zero. This example is recalled in the next Sections, in order to illustrate some definitions and theorems.

	As a final remark, note that definition (\ref{fvs}) does not give any hint on how to calculate the index $r$. Actually, the calculation of $r$ depends on the particular problem that we want to solve by using the concept of indexed identity. Any generalization of equation (\ref{calculandor}) does not necessarily encompass all possible methods for calculating $r$.

\section{`Fuzzy' Set Theory with Indexed Identity}

	In this Section we present an i-fuzzy set theory (`i' stands for indexed) based on the concept of indexed identity. We also show that this is a special case of Zadeh's original fuzzy set theory.

\begin{definicao}\label{fuzzyset}
If $X$ is an indexed variables system then an {\em i-fuzzy set\/} is a function $F:X\to R$, such that $\langle x,r\rangle\in F$ iff $\exists y\in X (y\equiv_r x)$.
\end{definicao}

	According to current literature \cite{Mordeson-98}, a fuzzy subset of a given $X$ is a function from $X$ into $[0,1]$. So:

\begin{teorema}
Every i-fuzzy set is a fuzzy set.
\end{teorema}

\noindent
{\bf Proof:} Straightforward from definitions of fuzzy set and i-fuzzy set.$\Box$\\

	This last theorem allows us to establish a relationship between fuzzy set theory (a la Zadeh) and our indexed variables system.

\begin{definicao}
The set of all i-fuzzy sets $F:X\to R$ is denoted by ${\cal F}(X;R)$.
\end{definicao}

\begin{definicao}
Let $X$ be an indexed variables system. If $F$ is an i-fuzzy set and $x$ is an element of $X$, then $x\in_r F$ iff $\langle x,r\rangle \in F$. In other words, $x\in_r F$ iff $F(x) = r$.
\end{definicao}

\begin{definicao}\label{union}
If $F$ and $G$ are i-fuzzy sets, the {\em union\/} $F\cup G$ is a function $F\cup G: X\to R$ defined as it follows: $\langle x,r\rangle\in F\cup G$ iff $r = \max \{F(x),G(x)\}$.
\end{definicao}

\begin{definicao}\label{intersection}
If $F$ and $G$ are i-fuzzy sets, the {\em intersection\/} $F\cap G$ is a function $F\cap G: X\to R$ defined as it follows: $\langle x,r\rangle\in F\cap G$ iff $r = \min \{F(x),G(x)\}$.
\end{definicao}

\begin{teorema}
The union of two i-fuzzy sets is an i-fuzzy set.
\end{teorema}

\noindent
{\bf Proof:} According to definition (\ref{union}), if $\langle x,r\rangle\in F\cup G$ then $r=F(x)$ or $r=G(x)$. Since $F$ is an i-fuzzy set, then there exists $y\in X$ ($X$ is a given indexed variables system) such that $y \equiv_{F(x)} x$. Analogously, there is $z\in X$ where $z \equiv_{G(x)} x$. So, there is $w$ (which is $y$ or $z$) such that $w\equiv_{\max \{F(x),G(x)\}}x$. Hence, $F\cup G$ is an i-fuzzy set.$\Box$

\begin{teorema}
The intersection of two i-fuzzy sets is an i-fuzzy set.
\end{teorema}

\noindent
{\bf Proof:} Analogous to the previous proof.$\Box$

\begin{definicao}\label{distinctsets}
Let $X$ be an indexed variables system. If $F$ and $G$ are i-fuzzy sets then the {\em distinction\/} between $F$ and $G$ is given by

$$D(F,G) = \sup_{x\in X}|F(x)-G(x)|,$$

\noindent
where $\sup$ stands for the supremum.
\end{definicao}

\begin{definicao}\label{equivsets}
If $F$ and $G$ are i-fuzzy sets then $F\equiv_r G$ iff $D(F,G) = 1-r$. The set of predicates $\equiv_r$ defined on elements of ${\cal F}(X;R)$ is denoted by $\Xi$.
\end{definicao}

\begin{teorema}
$\langle{\cal F}(X;R),\Xi\rangle$ is an indexed variables system.
\end{teorema}

\noindent
{\bf Proof:} We have to prove that $\langle{\cal F}(X;R), \Xi\rangle$ satisfies axioms {\bf F1}--{\bf F7}. So, we split this proof into seven parts. Hence, if $F$, $G$, and $H$ are i-fuzzy sets:

\begin{enumerate}

\item Since $X$ is nonempty, then ${\cal F}(X;R)$ is nonempty. This verifies axiom {\bf F1}.

\item According to definition (\ref{fuzzyset}), $0\leq F(x)\leq 1$ and $0\leq G(x)\leq 1$, for all $x\in X$. So, $0\leq\sup_{x\in X}|G(x)-F(x)|\leq 1$, i.e., $0\leq D(F,G)\leq 1$. Since $D(F,G) = 1-r$, where $F\equiv_r G$ (definition (\ref{equivsets})), then $\Xi = \{\equiv_r\}_{r\in R}$ is a family of binary predicates defined on the elements of ${\cal F}(X;R)$, where $R$ is a subset of the interval of real numbers $[0,1]$, such that $1\in R$. As a matter of fact, the proof that $1\in R$ is given in step 4 of this proof. This verifies axiom {\bf F2}.

\item $\sup_{x\in X}|F(x)-G(x)| = \sup_{x\in X}|G(x)-F(x)|$, i.e., $D(F,G) = D(G,F)$. So, according to definition (\ref{equivsets}), $F\equiv_r G$ iff $G\equiv_r F$. This verifies axiom {\bf F3}.

\item $F\equiv_1 G$ iff $F=G$, according to definitions (\ref{distinctsets}) and (\ref{equivsets}). This verifies axiom {\bf F4}.

\item If $F\equiv_r G$ then $D(F,G) = 1-r$. According to definition (\ref{distinctsets}) there is no $s$ such that $s\neq r$ and $D(F,G) = 1-s$. So, there is no $s$ such that $s\neq r$ and $F\equiv_s G$. This verifies axiom {\bf F5}.

\item Since $F$ and $G$ are limited functions, there is always a supremum $\sup_{x\in X}|G(x)-F(x)|$. So, $\forall x\forall y\exists r (F\equiv_r G)$. This verifies axiom {\bf F6}.

\item Since ${\cal F}(X;R)$ is a space of limited functions, then $\langle {\cal F}(X;R), D\rangle$ is a metric space, where $D$ is given as in definition (\ref{distinctsets}). This occurs because the distinction between two fuzzy sets is the well known {\em metric of uniform convergence\/} or {\em sup metric\/}. Then the triangle inequality is satisfied. This verifies axiom {\bf F7}.$\Box$

\end{enumerate}

\begin{corolario}
$D(F,G)$ is a distance between two i-fuzzy sets, and $\langle {\cal F}(X;R), \Xi\rangle$ induces a metric space $\langle {\cal F}(X;R), D\rangle$.
\end{corolario}

	Recalling the example given in Section III, we can define, e.g., the following i-fuzzy sets: $F = \{\langle 8,1\rangle, \langle 2,2/3\rangle, \langle 3,0\rangle\}$, $G = \{\langle 3,1\rangle, \langle 2,1/3\rangle, \langle 8,0\rangle\}$, and $H = \{\langle 2,1\rangle, \langle 8,2/3\rangle, \langle 3,1/3\rangle\}$. So, $F\cup G = \{\langle 8,1\rangle, \langle 3,1\rangle, \langle 2,2/3\rangle\}$, $F\cap G = \{\langle 3,0\rangle, \langle 8,0\rangle, \langle 2,1/3\rangle\}$, and $D(H,F\cup G) = 2/3$.

\section{Indexing Predicates}

	In this Section we give a general procedure which allows us to index any predicate, at least in principle.

	Since our mathematical framework is set theory, consider a set $S$ defined by means of the Separation Schema of Zermelo-Fraenkel set theory:

\[S = \{x\in X; P(x)\}\]

\noindent
where $X$ is a given set and $P$ is a given predicate. So, if $X$ is `the set of human beings' and $P$ is the predicate `to be smart', then $S$ corresponds to the set of smart people.

	The question now is: how to index predicate $P$? In other words: how to index the set $S$? The procedure that we suggest follows in the next paragraphs.

\begin{definicao}
Let $S$ and $X$ be the sets given above. If $x$ is an element of $X$, the distinction between $x$ and $S$ is
\[D(x,S) = \inf_{y\in S}D(x,y)\]
\end{definicao}

	Note that this last definition allows us to index the predicate $P$, i.e., the set $S$.

	{\em By indexing $P$ we mean the definition of $D(x,S)$.}

\begin{definicao}
$x\equiv_r S$ iff $D(x,S) = 1-r$.
\end{definicao}

\begin{exemplo}
If $X$ is `the set of human beings' and $P$ is the predicate `to be smart', then $S$ corresponds to the set of smart people. If we define a distinction function $D(x,S)$ where $x\in X$, then we are indexing the concept of `being smart'. If $D(Einstein, S) = 0.5$, then $Einstein \equiv_{0.5} S$, i.e., $Einstein$ is a half-smart person. The index $r$ corresponds to a degree of smartness.
\end{exemplo}

\begin{exemplo}
If $X$ is the set of subsets of a given set $T$ and $P$ is the predicate `to be open' (in the usual topological sense), then $S$ corresponds to the topology of a topological space. If we define a distinction function $D(x,S)$ where $x\in X$, then we are indexing the concept of `being open'. If $D(a, S) = 0.7$, then $a \equiv_{0.3} S$, i.e., $a$ has a 0.3 degree of openess.
\end{exemplo}

	These two examples can help us to see how powerful is this method of {\em indexation\/}. In the the second example we showed how to index the concept of `open set' in a topological space. But we can discuss about how to index the very concept of {\em topological space\/}.

	If we want to index the predicate `to be a topological space' rather than topological concepts like `to be open' or `to be compact', then we need: (1) a universe class $X$ which corresponds to the collection of ordered pairs $\langle T,{\cal T}\rangle$ of sets; and (2) a predicate $P$ such that $P(\Im)$ iff $\exists T\exists {\cal T}$ such that (i) $\Im = \langle T,{\cal T}\rangle$, (ii) $T$ is a non-empty set, (iii) the elements of ${\cal T}$ are subsets of $T$, (iv) $\emptyset\in{\cal T}$, (v) $T\in{\cal T}$, (vi) if $t_1$ and $t_2$ are elements of ${\cal T}$ then $t_1\cap t_2\in{\cal T}$, and (vii) an arbitrary union of elements of ${\cal T}$ is still an element of ${\cal T}$. Besides, we need a distinction function (which is a metric) $D:X\times X\to [0,1]$. Since in Zermelo-Fraenkel set theory there is no such a thing like the set of all ordered pairs of sets, then we cannot found our mathematical framework into usual set theory. We could consider $X$ as a {\em category\/}. In this case we should extend definition (\ref{fvs}) to a category-theoretical predicate, which is a task for future works.

	Something analogous could be said about groups, vector spaces, lattices, fields, and other mathematical theories usually founded into the scope of set theory.

\section{Conclusions}

	The main advantages of our mathematical framework are:

\begin{enumerate}

\item It is a kind of fuzzy mathematics (in the intuitive sense), which fuzzifies the concept of {\em equality\/} rather than that one of {\em membership\/} (as in the original work of Zadeh). So, it offers another point of view in the process of fuzzification. 

\item It allows us to use the theory of metric spaces, at least in principle, in order to derive theorems in fuzzy set theory.

\item It allows us a generalized method of `fuzzification' of predicates, in the sense given in the previous Section.

\end{enumerate}

\section{Acknowledgements}

	We akcnowledge with thanks Aur\'elio Sartorelli and Soraya R. T. Kudri for helpful suggestions and criticisms.

\end{document}